\documentclass{article}

\usepackage{amsmath}
\usepackage{amssymb}
\usepackage{latexsym}
\usepackage{theorem}
\usepackage{amscd}
\usepackage{epsfig}

\newtheorem{theorem}{Theorem}[section]			 
\newtheorem{lemma}[theorem]{Lemma}				 
\newtheorem{proposition}[theorem]{Proposition}	 
\newtheorem{corollary}[theorem]{Corollary}

{\theorembodyfont{\rmfamily}					 
\theoremstyle{plain}							 
\newtheorem{definition}[theorem]{Definition}	 
			 
\newtheorem{problem}[theorem]{Problem}			 
			 
}

\newcommand{\A}{\ensuremath{\mathcal A}}
\renewcommand{\L}{\ensuremath{\mathcal L}}

\newcommand{\I}{\ensuremath{\mathcal I}}
\newcommand{\E}{\ensuremath{\mathcal E}}

\newcommand{\C}{\ensuremath{\mathbb{C}}}

\newcommand{\Z}{\ensuremath{\mathbb{Z}}}

\newcommand{\R}{\ensuremath{\mathbb{R}}}

\newcommand{\lc}{\ensuremath{\ell c}}
\renewcommand{\to}{\ensuremath{\longrightarrow}}
\newcommand{\la}{\ensuremath{\lambda}}
\newcommand{\we}{\ensuremath{\wedge}} 
\renewcommand{\bar}[1]{\ensuremath{\overline{#1}}}
\newcommand{\rk}{\operatorname{rk}}
\newcommand{\nb}{\operatorname{\it nbb}}
\newcommand{\nc}{\operatorname{\it nbc}}
\newcommand{\poly}{\operatorname{\it poly}}
\newcommand{\supp}{\operatorname{supp}}
\newcommand{\mc}{\operatorname{mc}}

\newcommand{\nbc}{{\em nbc}}
\newcommand{\NBB}{{\em NBB}}
\newcommand{\qed}{\hfill \mbox{$\Box$}\medskip}
\newenvironment{ack}{\noindent {\bf Acknowledgements.}}{\par}
\newenvironment{proof}{\noindent {\it proof:}}{\qed \par}
\renewcommand{\theenumi}{\rm \theenumi}

\renewcommand{\theenumi}{\roman{enumi}}

\title{Combinatorial and algebraic structure in Orlik-Solomon 
algebras\thanks{revised July 31, 2000}} 

\author{Michael Falk} 
\date{}
\begin{document}
\maketitle

\begin{abstract}
The Orlik-Solomon algebra $\A(G)$ of a matroid $G$ is the free 
exterior algebra on 
the points, modulo the ideal generated by the circuit boundaries. On 
one hand, this algebra is a homotopy invariant of the complement of 
any 
complex hyperplane arrangement realizing $G$. On the other hand, some 
features of the matroid $G$ are reflected in the algebraic structure 
of $\A(G)$. 

In this mostly expository article, we describe recent developments in 
the construction of algebraic invariants of $\A(G)$. We develop a 
categorical framework for the statement and proof of recently 
discovered isomorphism theorems which suggests a possible setting for 
classification theorems. Several specific open problems are 
formulated. 
\end{abstract}

\begin{section}{Introduction: \\ The Orlik-Solomon algebra of a matroid}

Let $G$ be a simple matroid with ground set $[n]:=\{1,\ldots,n\}$. 
The Orlik-Solomon ($OS$) algebra of $G$ is defined as follows. Let 
$\E=\Lambda(e_1,\ldots,e_n)$ be the graded exterior algebra on 
elements 
$e_i$ of degree one corresponding to the points of $G$. For 
simplicity 
we will assume the ground field is $\C$. Except where noted, all of 
the results will hold for coefficients in an arbitrary commutative 
ring.

Define the linear mapping $\partial : \E^p \to \E^{p-1}$ by 
$$\partial(e_{i_1} \we \cdots
\we e_{i_p})=
\sum_{k=1}^p
(-1)^{k-1}e_{i_1} \we \cdots \we \widehat{e}_{i_k} \we \cdots \we 
e_{i_p},$$
where $\widehat{\hspace{1em}}$ indicates an omitted factor.

If $S=(i_1,\ldots,i_p)$ is an 
ordered $p$-tuple we denote the  product
$e_{i_1}\we\cdots \we e_{i_p}$ by $e_S$.
Let $\I$ denote the ideal of \E\ generated by $\{\partial e_S \ | \ S 
\ \text{is dependent}\}$. 
\begin{definition} The {\em Orlik-Solomon algebra} $\A=\A(G)$ of $G$ 
is 
the quotient $\E/\I$.
\end{definition}
Since \I\ is generated by homogeneous elements, both \I\ and \A\ 
inherit gradings from \E. We will denote the image of $e_S$ in $\A$ 
by $a_S$.

The $OS$ algebra has both combinatorial and topological significance, 
as demonstrated by these two results from \cite{OS1}. Recall that a 
projective realization of $G$ gives rise to a linear hyperplane 
arrangement. Throughout the paper $A$ will denote a hyperplane 
arrangement arising from a complex projective realization of $G$, 
and $M$ will denote the {\em complement} of $A$, $M=\C^\ell - 
\bigcup_{H\in A} H$.
\begin{theorem}The $OS$ algebra $\A(G)$ is isomorphic to the 
cohomology algebra $H^*(M)$.
\label{orsol}
\end{theorem}
The Whitney numbers of the second kind are defined in terms of the 
M\"obius function $\mu: L(G) \to \Z$ of the lattice of flats $L(G)$. 
Specifically, $$w_p(G)=\sum_{X\in L,\rk(X)=p}(-1)^p\mu(0_L,X).$$
\begin{theorem}The dimension of $\A^p(G)$ is equal to the $p^{\rm 
th}$ Whitney number $w_p(G)$ of $G$.
\label{whitney}
\end{theorem}

Theorem \ref{orsol} motivates what is for us the main problem 
concerning $OS$ algebras: to classify $\A(G)$ up to isomorphism of 
graded algebras. This type of problem is more familiar in topology 
than combinatorics, but the classification in this instance will be 
purely matroidal. Theorem~\ref{whitney} provides one line along which 
a classification could proceed, that is to extract combinatorial 
features of 
the matroid $G$ from algebraic invariants of $\A(G)$. In this regard 
we note that there are many sets of matroids with identical Whitney 
numbers 
while, on the other hand, the betti numbers $\dim(\A^p(G))$ in a 
sense 
take no account of the ring structure of $\A(G)$.

These observations set the tone for the exposition to follow.
We will construct multiplicative invariants of $\A(G)$ and attempt to 
extract combinatorial structure from them. The most delicate of these 
are the resonance varieties, discussed in Section \ref{resonant}. In 
Section~\ref{tutte} we show how ``stabilized'' parallel connection 
and 
direct sum of matroids yield isomorphic $OS$ algebras. 
We also show that truncations of matroids with isomorphic $OS$ 
algebras will have the same property. We make sense 
of these results using the categories of pointed matroids and affine 
$OS$ algebras, indicating a framework for the eventual 
classification. 
In Section \ref{lines} we describe recent work relating the $k$-adic 
closure of $\A(G)$ to the ``$k$-closure'' of the matroid $G$.

We close this introduction by recalling the oldest multiplicative 
invariant of $\A(G)$, termed ``the global invariant'' $\phi_3$ in 
\cite{F6}. Consider the multiplication map $$d: \E^1\otimes \I^2 \to 
\E^3.$$ This linear map can be shown to be an invariant of $\A(G)$. 
The nullity of $d$ is denoted by $\phi_3(\A)$. This quantity has a 
topological interpretation in terms of the 
fundamental 
group of the complement $M$. Indeed, the definition of $\phi_3$ comes 
directly out of the 
study of the rational homotopy type of hyperplane complements 
\cite{F4}. And of course $\phi_3(\A)$ can be thought of as an invariant 
of the 
matroid $G$. But the following problem remains open, even for graphic 
matroids.
\begin{problem}{Give a combinatorial interpretation of $\phi_3(G)$.}
\label{phi3}
\end{problem}
We will return to this problem in Section \ref{lines}.

The reader is referred to \cite{OT} for background material on 
complex 
hyperplane arrangements and Orlik-Solomon algebras, and to \cite{Wh1} 
for matroid theory. Section \ref{resonant} is largely based on
\cite{F7}, and much of Section \ref{tutte} is a reformulation of part 
of \cite{EF}. Section \ref{lines} is a brief report on work in 
progress; details and proofs will appear in \cite{F14} and \cite{DY}.
\end{section}

\begin{section}{Resonance varieties}
\label{resonant}

To answer questions concerning generalized hypergeometric functions, 
we began studying the $OS$ algebra as a differential complex in 
\cite{FT}, and then realized that our work could be used to define 
algebraic invariants \cite{F7}.

Fix an element $a_\la=\sum_{i=1}^n \la_i a_i$ 
in $\A^1$. Then left multiplication by $a_\la$ defines a map $\A^p 
\to \A^{p+1}$, which squares to zero. Thus we have a cochain complex
$$0\to\A^0 \overset{a_\la}{\to} \A^1 \overset{a_\la}{\to} \ldots 
\overset{a_\la}{\to} 
\A^{\ell-1} \overset{a_\la}{\to} \A^\ell\to 0.$$

The cohomology of this complex determines a stratification of the 
parameter space $\C^n$. The $p^{\rm th}$ {\em resonance variety} of 
\A\ is 
defined by 
$$R_p(\A)=\{ \la \in \C^n \ | \ H^p(\A,a_\la)\not = 0\}.$$ It is 
shown in \cite{F7} that $R_p(\A)$, up to ambient linear isomorphism, 
is an invariant of \A.

Basic properties of resonance varieties follow from 
the main results 
of \cite{Y2}. Let $\Delta$ denote the diagonal hyperplane 
$\sum_{i=1}^n \la_i = 0.$ Then

\begin{itemize}
\item $0\in R_p(\A)$ for $0\leq p\leq \ell$.
\item $R_0(\A)=\{0\}.$
\item $R_p(\A)\subseteq \Delta$ for all $p$.
\item $R_\ell(\A)\subseteq R_{\ell-1}(\A)$.
\item if $G$ is connected, then $R_\ell(\A)=R_{\ell-1}(\A)=\Delta$.
\item $R_p(\A)$ is a proper subvariety of $\Delta$ for $0\leq p\leq 
\ell-2$.
\end{itemize} 

Under some genericity conditions on \la, the cohomology 
$H^*(\A,a_\la)$ is isomorphic to the cohomology of $M$ with 
coefficients in a rank-one complex local system $\L_\la$ with monodromy 
determined by $\la$. This local system cohomology plays a role in the 
definition of generalized (multivariate) hypergeometric integrals. In 
a sense made precise in recent work of D.~Cohen and P.~Orlik 
\cite{CO}, the complex $(\A,a_\la)$ is the derivative at the identity 
of a cochain complex $(\A,\Delta_\la)$ that computes the local system 
cohomology. 
The resonance variety $R_p(A)$ is then the tangent cone at the identity to the 
``jumping locus'' for the local system cohomology, the set of local 
systems for which the cohomology $H^p(M,\L_\la)$ is non-vanishing. 
For $p=1$ the 
jumping locus for local system cohomology 
coincides with the character variety in $(\C^*)^n$ associated with the 
Alexander invariant of the fundamental group. For any $p$, a theorem of 
D.~Arapura asserts that these jumping loci are subtori of $(\C^*)^n$, possibly translated 
by elements of finite order.
This gives an indication of the proof of the following result, 
originally conjectured for $p=1$ in \cite{F7}, proved 
in that special case in \cite{CS3} and \cite{LY}, and finally 
established for arbitrary $p$ in \cite{CO} and \cite{L6}. See 
those papers for complete references. 

\begin{theorem} The resonance variety $R_p(\A)$ is a union of linear 
subspaces of $\C^n$.
\label{linear}
\end{theorem}

By Theorem \ref{linear}, $R_p(\A)$ can be thought of 
as a subspace arrangement, and as such, realizes a polymatroid 
$\poly_p(\A)$, which in essence records the dimension of the span of 
each subcollection of irreducible components of $R_p(\A)$. Because 
$R_p(\A)$ is invariant up to linear change of coordinates, the 
polymatroid $\poly_p(\A)$ is indeed an invariant of \A, powerful 
enough (at least for $p=1$) to distinguish $OS$ algebras of matroids 
which are almost identical in other respects \cite{F7}.

The first cohomology $H^1(\A,a_\la)$ can be computed directly, 
yielding a description of $R_1(\A)$. The following lemma reduces the 
calculation to an analysis of elements of $\I^2$.

\begin{lemma} $\la\in R_1(\A)$ if and only if $e_\la$ is one factor 
of a nonzero elementary tensor in $\I^2$.
\end{lemma}
Proof of this lemma and the results to follow can be found in 
\cite{F7}.

Irreducible components of $R_1(\A)$ are contained in 
intersections of 
$\Delta$ with hyperplanes $H_X$ defined by $\sum_{i\in X}\la_i=0$, 
where $X$ runs over certain flats of $G$. The flats which occur in 
these intersections are determined by so-called ``neighborly 
partitions'' of $G$.

\begin{definition} A {\em neighborly partition} of $G$ is a 
partition $\Pi$ of $[n]$ such that $|\pi \cap X|\not = |X|-1$ for all 
blocks $\pi\in\Pi$ and flats $X$ of rank two in $L$.
\end{definition}

We say a flat $X$ is 
``multi-colored'' if $X$ meets more than one block of $\Pi$. 
Given a neighborly partition $\Pi$ of a submatroid $S\subseteq [n]$ 
of $G$, set 
$$L_\Pi=\Delta \ \cap \ \bigcap_{i\not\in S}H_i \ \cap \ 
\bigcap_{X\in\mc(\Pi)} H_X,$$ where the 
last intersection runs over 
the set $\mc(\Pi)$  multi-colored rank-two flats of $\Pi$. Note that 
$H_i=\{\la\in 
\C^n \ | \ \la_i=0\}$.
The {\em support} $\supp(\la)$ of $\la$ is $\{i\in [n]  \ | \ 
\la_i\not=0\}$, considered as a submatroid of $G$. Let $\sim$ denote 
the equivalence relation associated with $\Pi$. Finally, for 
$\tau\in\E^2$ write $\tau =\sum_{i<j}\tau_{ij}e_i\we e_j$.
Here then is a description of $R_1(\A)$, from \cite{F7}, to which the 
reader is referred for the proof, examples and 
consequences.

\begin{theorem} $\la \in R_1(\A)$ if and only if $\supp(\la)$ affords 
a neighborly partition $\Pi$ such that (i) $\la \in L_\Pi$, and (ii) 
there exists $\mu\in L_\Pi$ not proportional to $\la$ 
such that $(\la\we\mu)_{ij}=0$ for every $i<j$ with
$i\sim j$ under $\Pi$.
\label{descrip}
\end{theorem}
The second condition will be replaced with a simpler criterion below.
 
If $X$ is a flat of rank two with $|X|\geq 3$, then 
$\Pi=\{\{i\} \ | \ i\in X\}$  is a neighborly 
partition of $X$, and $L_\Pi=\Delta \ \cap \ \bigcap_{i\not\in X} 
H_i$ has 
dimension $|X|-1\geq 2$. Thus condition (ii) is satisfied, and 
indeed $L_\Pi$ is a component of $R_1(\A)$ \cite{F7}. The 
components which 
arise in this way are called the {\em local components} of 
$R_1(\A)$. Here is a sample result from \cite{F7} showing how 
combinatorial structure may be extracted from $R_1(\A)$. 

\begin{corollary} Suppose every non-local component of 
$R_1(\A)$ has dimension two. Then $R_1(\A)$ determines the number of 
rank-two flats of $G$ of each cardinality. In particular, if $G$ has 
rank three, 
$R_1(\A)$ determines the Tutte polynomial of $G$.
\label{tuttegood}
\end{corollary}
D.~Cohen informs us that he and J.~Oxley have found examples for which the hypothesis fails.
We will see in the next section that  
$\A(G)$ does not generally determine the Tutte polynomial of $G$ for matroids of high rank.

\medskip
In \cite{LY} A.~Libgober and S.~Yuzvinsky base a study of the resonance variety 
$R_1(\A)$ on the Vinberg classification of Cartan matrices for 
affine Kac-Moody Lie algebras. Their approach yields substantial 
additional detail about $R_1(\A)$ and the associated 
neighborly partitions. We 
state some of their more general conclusions in the following theorem.

\begin{theorem}[\cite{LY}] \begin{enumerate}\item The irreducible components of 
$R_1(\A)$ are precisely the $L_\Pi$ of dimension at least two.
\item If $L_\Pi$ and $L_{\Pi'}$ are two irreducible components 
of $R_1(\A)$, then $L_\Pi\cap L_{\Pi'}=\{0\}$.
\item For any component $L_{\Pi}$ of $R_1(\A)$, each 
multi-colored flat of $G$ meets every block of $\Pi$.
\end{enumerate}
\label{yuz}
\end{theorem}
Theorem \ref{yuz}(i) effectively replaces condition (ii) of Theorem 
\ref{descrip} with the much simpler  requirement $\dim(L_\Pi)\geq 2$.
Theorem \ref{linear} for $p=1$ is an immediate 
corollary.

Matroids of rank greater than two
which support neighborly partitions $\Pi$ for which $L_\Pi$ has 
dimension at least two are quite rare. Some examples appear in 
\cite{F7}. The classification theory used in \cite{LY} imposes some 
restrictions, and also yields a method of constructing examples as a 
kind of inverse problem. The first part of the following problem is solved in 
some special cases in \cite{LY}.

\begin{problem}\begin{enumerate} \item Characterize those matroids 
which support 
neighborly partitions $\Pi$ satisfying $\dim(L_\Pi)\geq 2$.
\item Describe the polymatroid $\poly_1(G)$ associated with the 
arrangement of subspaces $\{L_\Pi \ | \ \Pi \ \text{is neighborly 
and} \ \dim(L_\pi)\geq 2\}$.
\end{enumerate}
\label{neighborly}
\end{problem}

Libgober and Yuzvinsky \cite{LY} also uncover 
a connection between
non-local components of $R_1(\A)$, for arrangements of rank three, and 
pencils of curves $\C P^2 \to \C P^1$ which include the arrangement 
in their singular locus. The existence of such pencils imposes 
further restrictions on the structure of matroids supporting 
nontrivial ($\dim(L_\pi)\geq 2$) neighborly partitions. In addition, 
these pencils of curves
bear a relationship to the $K(\pi,1)$ problem for complex hyperplane 
arrangements, and were studied in that vein in \cite{F10}. So a 
solution to Problem \ref{neighborly}(i) might have some implications 
for the $K(\pi,1)$ problem \cite{FR4}.

\medskip
In another direction, D.~Matei and A.~Suciu \cite{MatSuc} discovered deep connections 
between the resonance varieties of $\A(G)\otimes \Z_p$ and 
the structure of the second nilpotent quotient of the fundamental 
group $\pi_1(M)$. 
This work leads to some other interesting open 
questions. We briefly summarize.

Write $R_1(\A,\Z_p)$ for the first resonance variety of 
$\A(G)\otimes \Z_p$, and let $$R_{1,d}(\A,\Z_p)=\{\la \in 
R_1(\A,\Z_p) \ | \ \dim H^1(\A\otimes \Z_p, a_\la)\geq d\}.$$
These are subvarieties of $(\Z_p)^n$, easily seen to be homogeneous. 
Let $\widehat{R}_{1,d}(\A,\Z_p)$ denote the projective image of 
$R_{1,d}(\A,\Z_p)$. Finally, let $\pi=\pi^1\supseteq \pi^2 
\supseteq\pi^3\supseteq \cdots$ denote the descending central series 
of $\pi=\pi_1(M)$, and $\Gamma=\pi/\pi^3$ the second 
nilpotent quotient. Let $\nu_{p,d}$ denote the number of normal subgroups $K$ of 
$\Gamma$ of index $p$, such that the abelianization of $\Gamma/K$ 
has $p$-torsion of rank $d$.

\begin{theorem}[\cite{MatSuc}] 
$$\nu_{p,d}=\sharp\left(\widehat{R}_{1,d}(\A,\Z_p) - 
\widehat{R}_{1,{d+1}}(\A,\Z_p)\right)$$
\label{matsuc}
\end{theorem}
The quantity on the right-hand side is also an invariant of $\A(G)$. 

The proof uses a relationship between 
the resonance varieties and the Alexander invariant of the fundamental 
group, similar to the observations used to prove Theorem \ref{linear} 
in \cite{CS3}. In this case, the (linearized) Alexander matrix 
(mod $p$) is used to count normal subgroups of index $p$ in the 
second nilpotent quotient of $\pi_1(M)$, on one hand, and to define 
the resonance variety of $\A(G)\otimes \Z_p$ on the other. 

Theorem \ref{matsuc} leads to the study of resonance varieties of 
$OS$ algebras over 
finite fields. Because the variety $R_1(\A)$ is defined over \Z, we 
can reduce mod $p$. But there are matroids $G$ which have ``exceptional primes'' 
$p$, for which the reduction
$R_1(\A)\otimes \Z_p$ does not coincide with $R_1(\A,\Z_p)$. The 
basic results of this section, from \cite{F7}, will hold over an 
arbitrary ground field, but the techniques of \cite{LY} and 
\cite{CS3}, for instance, and thus Theorems \ref{linear} and \ref{yuz}, 
require complex coefficients. In 
\cite{MatSuc} the authors give examples of matroids for which
\begin{enumerate} 
\item $R_1(\A,\Z_p)$ has non-local components while $R_1(\A)$ has 
none. 
\item $R_1(\A,\Z_p)$ has a non-local components of dimension greater 
than two, while all non-local components of $R_1(\A)$ are $2$-dimensional.
\item $R_{1,d}(\A,\Z_p)$ has components which are not $(d+1)$-dimensional. 
By contrast, the components of the 
analogous variety $R_{1,d}(\A)$ over \C\ always have dimension $d+1$ 
\cite{LY}.
\end{enumerate}

This suggests a variation of Problem \ref{neighborly}, 
suggested by A.~Suciu.
\begin{problem} Given a matroid $G$, determine the exceptional primes 
for $G$, that is, the primes $p$ for which $\R_1(\A,\Z_p)\not \cong 
\R_1(\A)\otimes \Z_p$.
\end{problem}
\end{section}

\begin{section}{Isomorphisms: Affine $OS$ algebras and pointed 
matroids}
\label{tutte}

In \cite{EF} we showed how one could construct, from an arbitrary 
pair 
of (realizable) matroids $G_0$ and $G_1$, a pair of non-isomorphic 
matroids $G$ 
and $G\,'$ for which $\A(G)\cong\A(G\,')$. The matroids $G$ and 
$G\,'$ 
are, respectively, the direct sum $G_0\oplus G_1$, and any parallel 
connection $P(G_0,G_1)$, stabilized  by adding an isthmus (so $G$ and 
$G\,'$ have the 
same number of points). In this section we cast this 
result in a simpler conceptual framework, motivated by the fact that 
parallel connection is the categorical direct sum of base-pointed 
matroids \cite{Bry71,Wh1}. 

We will also prove that, for two matroids $G$ and $G\,'$, if 
$\A(G)\cong\A(G\,')$, then $\A(\bar{G})\cong \A(\bar{G\,'})$, where 
the bar denotes truncation. Together with the 
equivalences 
involving direct sum, this result explains all known instances of 
isomorphisms of $OS$ algebras, and so we are led to a 
possible formulation for a classification result.

\medskip
We start with some fundamental observations. The elementary proofs 
are left to	the reader.									
											
\begin{proposition}
{\begin{enumerate}
\item If $i\in S$ then $e_i\partial e_S=\pm e_S$.
\item If $S$ is	dependent then $e_S\in \I.$					
\item The ideal	\I\ 
is generated by $\{\partial e_C \ | \ C \ \text{is a 
circuit}\}$.												  
\end{enumerate}}									
\label{elem}
\end{proposition}
									
\medskip
Our setup involves generalizing the definition of $OS$ algebra. This 
is carried out in \cite{OT} by giving an algebra presentation 
associated with an arrangement of affine hyperplanes. We adopt a 
different approach, so that we can stay in the realm of matroid 
theory. The combinatorial model for an affine arrangement is a 
{\em pointed matroid}, that is, a matroid with a specified base 
point. 
Given an arrangement $A$ of affine hyperplanes, the underlying 
pointed 
matroid will be the matroid of the cone $cA$ of $A$ \cite{OT}, with 
the hyperplane at infinity as base point. Conversely, given a central 
arrangement $A$ realizing the matroid $G$, the effect of choosing a 
base point in $G$ will yield the pointed matroid associated with the 
decone $dA$ of $A$ relative to the hyperplane corresponding to the 
chosen base point. In keeping with the notation of \cite{OT}, we will 
write $dG$ to denote a pointed matroid, with underlying unpointed 
matroid $G$. Our convention will be that $G$ has ground set 
$\{0,\ldots,n\}$, and that $dG$ has $0$ as base point. More 
generally, the pointed matroid on $G$ with base point $i$ will be 
denoted $d_iG$.

\begin{definition} The Orlik-Solomon ($OS$) algebra of the pointed 
matroid $dG$ is the subalgebra $\A_d(dG)$ of the $OS$ algebra 
$\A(G)$ generated by $\{a_1-a_0,\ldots,a_n-a_0\}$.
\label{affine}
\end{definition}

The reader will find that this definition agrees with the definition 
of \cite{OT} of the $OS$ algebra of an affine arrangement $dA$ with 
underlying pointed matroid $dG$. In particular we have \cite[Corollary 
3.58]{OT} $$\sum_p \dim(\A^p(G))t^p=(1+t)\sum_p \dim(\A_d(dG))t^p.$$

We recover the ordinary $OS$ algebra as follows. Given an unpointed 
matroid $G$ on ground set $[n]$, let $cG$ denote the matroid 
$\{0\}\oplus G$ of rank $\rk(G)+1$, with the point $0$ marked. Here 
$\{0\}$ is understood to be the rank-one matroid with one point, an 
isthmus. The reader is invited to 
verify the following result.

\begin{lemma} $\A_d(cG)\cong\A(G)$.
\label{cone}
\end{lemma}

There are two operations on pointed matroids which have a 
predictable effect on $OS$ algebras. The first of these will be 
obvious to those familiar with the topology of hyperplane 
arrangements. 
Indeed, the complement $M$ supports an action of $\C^*$, and the 
induced map $H^*(M/\C^*)\to H^*(M)$ is a split injection with image 
$\A_d(d_iG)$, for any $i$ 
\cite[Prop. 5.1]{OT}. 

\begin{theorem} For any $i,j\in\{0,\ldots,n\}$, 
$$\A_d(d_iG)=\A_d(d_jG).$$
\label{point}
\end{theorem}

\begin{proof} This is immediate from the identities
 $a_k-a_j=(a_k-a_i)-(a_j-a_i)$ for $k\not=i,j$ 
and $a_i-a_j=-(a_j-a_i)$.
\end{proof}

The parallel connection of pointed matroids $dG_0$ and $dG_1$ is the 
unique (up to isomorphism)
pointed matroid $P_d(dG_0,dG_1)$ of largest rank which is a union of 
pointed submatroids isomorphic to $dG_0$ and $dG_1$, whose ground 
sets intersect only at the 
base point \cite{Wh1}. The underlying 
matroid of $P_d(dG_0,dG_1)$ is called a parallel connection of $G_0$ 
and $G_1$, denoted $P(G_0,G_1)$. The following result from 
\cite{Bry71} motivated the present formulation of the equivalence 
discovered in \cite{EF}.

\begin{lemma} Parallel connection is a sum in the category of 
pointed matroids and pointed strong maps. That is,
\begin{equation*}
\begin{CD}
\{0\}@>>>G_0\\
@VVV @VVV\\
G_1@>>>P(G_0,G_1)
\end{CD}
\end{equation*}
is a pushout diagram of pointed strong maps.
\end{lemma}

\begin{lemma} The assignment $dG \mapsto \A_d(dG)$ yields a functor 
from the category of pointed matroids and pointed strong maps to the 
category of connected (i.e., $\A^0\cong \C$) graded algebras over \C.
\label{funct}
\end{lemma}

\begin{proof} Let $dG$ and $dG'$ be pointed matroids on 
$\{0,\ldots,n\}$ and $\{0,\ldots,m\}$  respectively.
A pointed strong map $dG \to dG'$ arises from a set 
function $\eta: \{0,\ldots,n\} \to \{0,\ldots,m\}$ mapping $0$ to 
$0$. 
This function yields a homomorphism of exterior algebras 
$\hat{\eta}: \E \to \E'$ determined by $\hat{\eta}(e_i)=e_{\eta(i)}$. 
According to \cite[Lemmas 8.1.4 and 8.1.6]{Wh1}, the image of each 
circuit of $G$ is dependent in $G'$. Using Lemma \ref{elem} this 
implies that $\hat{\eta}$ sends $\I$ 
into $\I'$, inducing a homomorphism $\A(G) \to \A(G')$. Since 
$\hat{\eta}(a_0)=a_0'$, $\hat{\eta}$ restricts to a homomorphism 
$\A_d(dG)\to \A_d(dG')$.
\end{proof}

As a consequence of these observations, the effect of parallel 
connection 
on $OS$ algebras becomes natural. 

\begin{theorem} The $OS$ algebra of $P_d(dG_0,dG_1)$ is isomorphic to 
$\A_d(dG_0)\otimes\A_d(dG_1)$.
\label{sum}
\end{theorem}

\begin{proof} Let us write $dG$ for $P_d(dG_0,dG_1)$. Using the fact 
that tensor product is a sum in the category of connected 
graded algebras, together with Lemma \ref{funct}, we obtain a 
surjective homomorphism $\A_d(dG_0)\otimes\A_d(dG_1) \to \A_d(G)$. 
Using Theorem~\ref{whitney} and \cite[Prop. 7.2.9]{Wh3}, one can show that 
the domain and target have the same dimension in each degree. Thus the 
two algebras are isomorphic.
\end{proof}

As a consequence of Theorems \ref{sum} and \ref{point}, we easily 
obtain the combinatorial/algebraic version of the main topological 
result of \cite{EF}.

\begin{theorem} Let $G_0$ and $G_1$ be arbitrary matroids. Then 
$G=G_0\oplus G_1$ and $G\,'=\{0\} \oplus P(G_0,G_1)$ have isomorphic 
$OS$ algebras.
\label{parallel}
\end{theorem}

\begin{proof} Consider the pointed parallel connection 
$d\hat{G}=P_d(cG_0,cG_1)$. The underlying matroid $\hat{G}$ is 
$\{0\}\oplus G_0\oplus G_1=\{0\}\oplus G$, which is precisely $cG$. 
Then, by  Lemma \ref{cone}, the $\A_d(dG)\cong\A(G)$. On the other 
hand, by Theorem \ref{sum}, $\A_d(dG)$ is also isomorphic to 
$\A_d(cG_0)\otimes \A_d(cG_1)$, which again by Lemma \ref{cone}, is 
isomorphic to $\A(G_0)\otimes \A(G_1)$.

Now, according to Theorem \ref{point}, we may change the base points 
of $cG_0$ and $cG_1$ without affecting the affine $OS$ algebras. The 
pointed parallel connection $d\hat{G}\,'$ of these new pointed 
matroids will have underlying matroid $\hat{G}\,'$ isomorphic to the 
the sum of two isthmuses (neither marked) with an ordinary parallel 
connection $P(G_0,G_1)$ of $G_0$ and $G_1$ along the new marked 
points of each. Again, we have $\A_d(d\hat{G}\,')\cong 
\A_d(cG_0)\otimes \A_d(cG_1)\cong \A(G_0)\otimes\A(G_1)$. Now we 
change the base point of $d\hat{G}\,'$ to one of the isthmuses, and 
recognize the resulting pointed matroid as $cG\,'$. We apply Lemma 
\ref{cone} once more to obtain the result.
\end{proof}

We regard the method of proof above as ``diagrammatic," and indeed 
the argument is easier to follow in pictures than in words. See 
Figure \ref{diagrams}.

%\psdraft
\begin{figure}[h]
\begin{center}
\epsfig{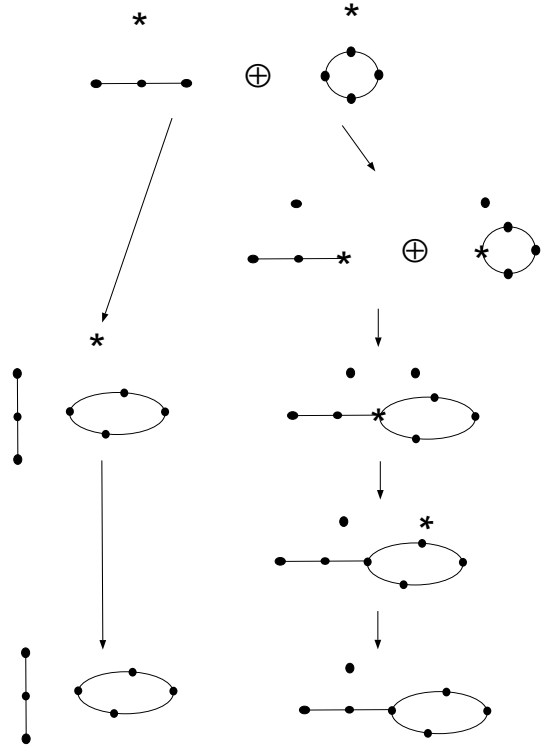}
\caption{The proof of Theorem \ref{parallel}}
\label{diagrams}
\end{center}
\end{figure}
%\psfull

It should now be clear that these isomorphisms arise from the
trivial operations of changing base points and forming sums.

In \cite{EF} we proved a stronger result for realizations of $G_0$ 
and $G_1$, by constructing a natural realization of $P(G_0,G_1)$ and 
proving that the complements of the arrangements realizing $G$ and 
$G\,'$ are in fact diffeomorphic. Theorem \ref{parallel} follows in 
this case by Theorem \ref{orsol}.

We state two interesting consequences of Theorem \ref{parallel} from 
\cite{EF}. The first should be compared with Theorems \ref{whitney} 
and \ref{tuttegood}.

\begin{corollary} Given an arbitrary matroid $G_0$, there exist 
extensions $G$ and $G\,'$ of $G_0$ with isomorphic $OS$ algebras but 
different Tutte polynomials.
\label{tuttediff}
\end{corollary}

The second corollary results from the indeterminacy in the change of 
base point in the proof of Theorem \ref{parallel}.

\begin{corollary} For any positive integer $n$, there exist $n$ 
nonisomorphic matroids with isomorphic $OS$ algebras.
\end{corollary}

The original examples of nonisomorphic matroids with isomorphic $OS$ 
algebras, which appeared in \cite{F6,OT,F3}, are truncations of $G$ 
and $G\,'$, where the factors $G_0$ and $G_1$ both have rank two. In 
an NSF-sponsored {\em REU} undergraduate research project directed by 
the author, C.~Pendergrass showed that truncation of matroids always 
preserves isomorphisms of the associated $OS$ algebra \cite{Pender}. 

\begin{theorem} Suppose $\A(G)\cong\A(G\,')$, and let $\bar{G}$ and 
$\bar{G\,'}$ denote the (corank-one) truncations of $G$ and $G\,'$ 
respectively. Then $\A(\bar{G})\cong \A(\bar{G\,'})$.
\label{trunc}
\end{theorem}

\begin{proof} Suppose $\eta$ is an isomorphism of $\A(G)$ to $\A(G\,')$. 
To begin with, we can then assume without loss that 
$G$ and $G\,'$ have the same ground set. 
The isomorphism $\eta: \A^1(G)\to 
\A^1(G\,')$ determines an 
isomorphism $\hat{\eta}: \E(G) \to \E(G\,')$, and 
$\hat{\eta}(\I(G)) = \I(G\,').$ We need only show that 
$\hat{\eta}(\I(\bar{G})) = \I(\bar{G\,'}).$

Let $n=\rk(G)=\rk(G\,')$. Then, for $p<n-1$,
$$\hat{\eta}(\I^p(\bar{G})) =\hat{\eta}(\I^p(G))
 =\I^p(G\,')= \I^p(\bar{G\,'}).$$  
Since the truncations have rank $n-1$, we also have, for $p\geq n-1$,
$\I^p(\bar{G})=\partial \E^{p+1}=\I^p(\bar{G\,'})$. Since 
$\hat{\eta}$ is an algebra homomorphism, it commutes with $\partial$, 
and thus $\hat{\eta}(\I^p(\bar{G})) = \I^p(\bar{G\,'})$ for $p\geq 
n-1$. This 
completes the proof.
\end{proof}

All known examples of isomorphisms of $OS$ algebras arising from 
nonisomorphic matroids are consequences of Theorems \ref{parallel} 
and \ref{trunc}. So we are led to the following problem. Recall that a 
matroid which is not a truncation is called {\em inerectible}.

\begin{problem} For inerectible parallel-irreducible 
matroids $G$ and $G'$, $\A(G)\cong\A(G\,')$ if and only if $G\cong G\,'$.
\end{problem}

We prefer an alternate formulation based on the categorical framework 
developed earlier.

\begin{problem} Suppose $dG$ and $dG'$ are inerectible pointed 
matroids which are irreducible in the category of pointed matroids. 
Then $\A_d(dG)\cong\A_d(dG')$ if and only if $dG\cong dG'$ up to 
change of base point.
\end{problem}

\end{section}
\begin{section}{The $k$-adic closure of $\A(G)$}
\label{lines}

We have recently become interested in quadratic $OS$ algebras, and 
more generally the quadratic closure of $\A=\A(G)$. This is the first in 
a series of {\em $k$-adic closures} whose dimensions are algebraic 
invariants of \A, and about which little is known. In this 
section we briefly present these ideas and 
describe some recent results and work in progress, to appear in 
\cite{F14} and \cite{DY}.

For $k\geq 2$, define the {\em $k$-adic $OS$ ideal} $\I_k$ to be the 
ideal generated by $\sum_{j\leq k} \I^j$ and the {\em $k$-adic closure 
of \A} to be the quotient $\A_k=\E/\I_k$. These algebras form a sort of 
resolution of \A:

$$\E=\A_1 \to \A_2 \to \A_3 \to \cdots \to \A_{\ell-1} \to 
\A_\ell=\A.$$

The following problem is wide open, even for $k=2$.

\begin{problem} Calculate the dimension of $\A_k^p$  in terms of the 
underlying matroid $G$.
\label{dim}
\end{problem}

Of special interest is the condition $\A_2=\A$, in which case we 
say $\A$ is {\em quadratic}. Examples indicate that this condition is 
related to the notion of line-closed matroid. The {\em line-closure} 
of a set $S\subseteq [n]$ is the smallest subset $\lc(S)$ of $n$ containing 
$S$ and containing the entire line in $G$ spanned by any pair of 
points of $\lc(S)$. The matroid $G$ is line-closed if and only if 
every line-closed set is closed. A proof of the following result 
will appear in \cite{F14}.

\begin{theorem} If \A\ is quadratic then $G$ is line-closed.
\label{quad}
\end{theorem}

This result was originally announced in \cite{F13}, at which time we 
conjectured that the converse is also true, that is, that line-closed 
matroids have quadratic $OS$ algebras. S. Yuzvinsky subsequently
found a counterexample to this conjecture, the matroid 
on eight points with nontrivial lines $$123, \ 3456, \ 167, \ 258, \ 
\text{and} \ 
478.$$ Yuzvinsky proposed a different 
condition for quadraticity of \A, which fails for the example 
above. This condition is also necessary for 
quadraticity, and is demonstrably stronger than line-closure. 
G.~Denham subsequently found an 
example (a $9_3$ configuration) showing this stronger condition is 
still not 
sufficient for 
quadraticity. The work of Denham and 
Yuzvinsky is	based on a detailed	study of the annihilator of	the	
quadratic $OS$ ideal $I_2$ inside the full tensor algebra, and is 
reported on in \cite{DY}. At this point there seems to be no easily 
stated matroidal criterion	equivalent to quadraticity.

\medskip
Theorem \ref{quad} is actually a corollary of a more general result 
concerning $\A_2$. We define a set $\nb(G)$ of increasing subsets of $[n]$ 
by $S=(i_1,\ldots,i_p)_<\in \nb(G)$ if and only if $i_j=\min 
\lc(\{i_j,\ldots, i_p\})$ for all $j$. This is an analogue of the set $\nc(G)$ 
of \nbc\ (=``no-broken-circuit'') sets of $G$ \cite{Bj2}. In fact these sets are 
precisely the \NBB\ (=``not-bounded-below'') sets of A.~Blass and B.~Sagan 
\cite{BSag}, which generalize \nbc\ sets, for the lattice of line-closed 
sets of $G$, with a linear ordering of the atoms. It is the case that 
$\nb(G)=\nc(G)$ if and only if $G$ is line-closed. Then \ref{quad} 
follows easily from the next theorem.

\begin{theorem} The set of monomials $\{a_S \ | \ S \in \nb(G)\}$ 
forms 
a linearly independent subset of $\A_2$.
\label{indep}
\end{theorem}

This generalizes half of the well-known theorem \cite{JT2,Bj2} that $\{a_S \ | \ 
S\in \nc(G)\}$ yields a basis for the $OS$ algebra \A.
Yuzvinsky's example shows that the set $\{a_S \ | \ S \in \nb(G)\}$
cannot form a basis for $\A_2$ in general. 

An analogue of Theorem \ref{indep} holds 
for $\A_k$ for each $k\geq 2$, giving a partial solution to Problem 
\ref{dim} in the form of combinatorial lower bounds. Of course, a 
formula for the cardinality of $\nb(G)$ has not been found. In fact, 
this cardinality can change if the linear order of the points is 
changed.

\begin{problem} Calculate the maximal cardinality of $\nb(G)$ over 
all linear orderings of the points of $G$.
\end{problem}

L.~Paris has informed us that 
$\{\partial e_C \ | \ C \ \text{is a circuit}\}$ can be shown directly to be a 
Gr\"obner basis for the $OS$ ideal \I. A complete direct proof is seemingly not 
extant. 
The fact that \nbc\ monomials form a basis for  \A\ is an immediate consequence. 
In fact the latter assertion implies the former -- see \cite[Theorem 4.1]{Pe}. 
The following problem seems more delicate.

\begin{problem} Find a Gr\"obner basis for the quadratic $OS$ ideal 
$\I_2$.
\end{problem}

Our experiments lead us to another interesting question, which seems 
to be related.

\begin{problem} Determine conditions on $S$ under which $\partial e_S$ 
will lie 
in the $k$-adic $OS$ ideal $\I_k$.
\end{problem}

We close by returning to the invariant $\phi_3$ defined in the 
Introduction. It turns out that a calculation of 
$\dim(\A_2^3)$ would yield a combinatorial formula for $\phi_3$.
Indeed, $\phi_3$ is the nullity of $\E^1\otimes \I^2 \to \E^3$, while 
the cokernel of the same map is precisely $\A_2^3$. The dimension 
of $\I^2$ is just $\dim(\E^2)-\dim(\A^2)$, so we obtain the following 
formula.

\begin{theorem} Let $n=\rk(G)$ and $w_2=\dim(\A^2)$, the second 
Whitney number of $G$. Then 
$$\phi_3= 2 {{n+1}\choose {3}} - nw_2 + \dim(\A_2^3).$$
\end{theorem}

Thus Problem \ref{phi3} is a special case of Problem \ref{dim}.

\end{section}

\medskip
\begin{ack} I am grateful to Raul Cordovil and Michel Las Vergnas for 
inviting me to speak at the CIRM conference, and for providing 
financial support. I thank Sergey Yuzvinsky for his help in studying 
quadratic $OS$ algebras, and Alex Suciu for his suggestion to include 
resonance varieties over $\Z_p$, and his help in understanding them. 
My {\em REU} students Carrie 
Eschenbrenner, Cayley Pendergrass, and Samantha Melcher assisted me
in sorting out much of the material in the last two sections. Finally I wish to 
thank Diane MacLagan for a helpful correspondence concerning Gr\"obner 
bases.
\end{ack}

%\bibliography{biblio}
%\bibliographystyle{plain}

\makeatletter
\bigskip \bigskip
{\small
Department of Mathematics and Statistics \\
Northern Arizona University \\
Flagstaff, AZ 86011-5717 \\
michael.falk@nau.edu 
}
\makeatother
 
\end{document}